\documentclass[manmat,numbook,envcountsame]{svjour}
\usepackage{amsmath,latexsym,amssymb,times,mathptm}

\def\sqr#1#2{{\vcenter{\hrule height.#2pt
    \hbox{\vrule width.#2pt height#1pt \kern#1pt
        \vrule width.#2pt}
    \hrule height.#2pt}}}
\def\square{\mathchoice\sqr64\sqr64\sqr{4}3\sqr{3}3}
\def\demo{\begin{proof}}
\def\QED{\hfill$\square$ \end{proof}}
\def\tratto{\mbox{\rule{2mm}{.2mm}$\;\!$}}
\def\lar{\longrightarrow}
\def\rar{\rightarrow}
\def\p{{\frak p}}

\def\m{{\frak m}}
\def\QQ{{\frak Q}}
\def\pp{{\frak p}}

\def\N{{\cal N}}
\def\L{{\cal L}}

\begin{document}

\title{On the integral closure of ideals\thanks{Part of the results
contained in this paper were obtained while the first author was
visiting Rutgers University and was partially supported by CNR
grant 203.01.63, Italy. The second and third authors were
partially supported by the NSF.} }

\titlerunning{On the integral closure of ideals}
\authorrunning{A. Corso et al.}

\author{Alberto Corso \and Craig Huneke \and Wolmer
V. Vasconcelos}

\mail{A. Corso, corso@math.purdue.edu.}

\institute{{\sc A. Corso} \\
Department of Mathematics, Purdue University, West Lafayette, IN
47907, USA \\ (e-mail$\colon$ corso@math.purdue.edu) \and
{\sc C. Huneke} \\
Department of Mathematics, Purdue University, West Lafayette, IN
47907, USA \\(e-mail$\colon$ huneke@math.purdue.edu) \and
{\sc W.V. Vasconcelos} \\
Department of Mathematics, Rutgers University, New Brunswick, NJ 08903, USA \\
(e-mail$\colon$ vasconce@rings.rutgers.edu)}

\date{Received July 28, 1997}

\maketitle

\begin{abstract}
Among the several types of closures of an ideal $I$ that have been
defined and studied in the past decades, the integral closure
$\overline{I}$ has a central place being one of the earliest and most
relevant. Despite this role, it is often a difficult challenge to
describe it concretely once the generators of $I$ are known.
Our aim in this note is to show that in a broad class of ideals their
radicals play a fundamental role in testing for integral closedness, and
in case $I\neq \overline{I}$, $\surd{I}$ is still helpful in finding
some fresh new elements in $\overline{I}\setminus I$.
Among the classes of ideals under consideration are: complete
intersection ideals of codimension two, generic complete
intersection ideals, and generically Gorenstein ideals.
\end{abstract}

\section{Introduction}
Let $R$ denote a Noetherian ring and $I$ one of its ideals.
The {\em integral closure} of $I$ is the ideal $\overline{I}$ of all
elements of $R$ that satisfy an equation of the form
\[
X^n+a_1X^{n-1}+\cdots+a_{n-1}X+a_{n}=0, \qquad a_i\in I^i.
\]
The literature is very spare on methods to find this ideal from
$I$, except when $R$ is a ring of polynomials over a field and $I$
is an ideal generated by monomials: $\overline{I}$ is then the
monomial ideal defined by the integral convex hull of the exponent
vectors of $I$ $($see \cite[p. 140]{Eisenbudbook}$)$. Already here
the problem has computationally a non-elementary solution. Thus
the general problem can be likened to solving a nonlinear integer
programming question---quite a tall order.

By comparison, the task of computing the {\em radical}, $\surd{I}$, of
the ideal $I$ is simpler as one is only required to describe the
solutions in $R$ of equations of the form
\[
X^m-b=0, \qquad b\in I.
\]

We always have
\[
I \subseteq \overline{I} \subseteq \surd{I},
\]
but in general these ideals differ from one another. A twin
problem, but much simpler to approach, consists in deciding the
equality $I = \overline{I}$. Such ideals are called {\it
integrally closed} or {\it complete}.

There are two basic ways to approach these questions.
A theoretical way of testing whether or not an element $z$ belongs to
the integral closure of an ideal is based on the so called {\em
determinant trick}: an element $z\in R$ is integral over $I$ if and
only if there exists a finitely generated faithful $R$-module $M$
such that $zM \subseteq IM$. In fact, for any such $M$ one has that
\[
IM \colon M \subseteq \overline{I}.
\]
In particular, if the ideal $I$ is integrally closed then
\begin{equation}\label{test}
\mbox{$IM \colon M=I$,}
\end{equation}
for any finitely generated faithful $R$-module $M$. The issue is to
find appropriate modules for a given ideal $I$. Some of the natural
choices we are going to employ are the powers of the ideal and of its
radical, Jacobian ideals, and modules of syzygies.

Another crude theoretical approach is through the theory of the Rees
algebra of $I$
\[
{\mathcal R} = R[It]= R + It+ I^2t^2 + \cdots \subset R[t],
\]
and one looks for its integral closure in $R[t]$,
\[
R + \overline{I}t + \overline{I^2}t^2 + \cdots \subset R[t].
\]
This is obviously wasteful of resources since the integral closure of
all the powers of $I$ will be computed.
Nevertheless, there are some instances when this path can be
profitably taken. First, in some special cases, it turns out that
$\overline{I^n}= (\overline{I})^n$ for all $n\geq 1$. Then, notions
typical of the theory of Rees algebras can be brought in. For example,
we say that $J \subseteq I$ is a
{\em reduction} of $I$ if $I^{r+1}=JI^r$ for some non-negative integer
$r$. A reduction is said to be {\em minimal} if it is minimal with respect
to containment. It is known that an ideal has the same integral
closure as any of its reductions $($see \cite{NR}$)$. Appropriately we can
then view all these algebras as modules over the Rees algebra $R[Jt]$.

\medskip

Our goals can be framed in terms of the following two issues:

\begin{description}[()]
\item[$\bullet$]
development of effective integral closedness criteria, i.e., tests
for the equality $I=\overline{I}$;

\item[$\bullet$]
development of the means to find the integral closure when the ideal $I$
fails to pass the above tests.
\end{description}

Our results do not meet this ambitious program, except for some
classes of ideals. In Sections $2$ and $3$, we give a very direct validation
of the equality $I=\overline{I}$ in two cases of interest:
$(${\it i}$)$ $I$ is a generic complete intersection, or $(${\it
ii}$)$ $I$ is a generically Gorenstein ideal. In these cases we provide a
`closed formula' $($see Theorem~\ref{ccitest1} and
Theorem~\ref{new}$)$ which says that $I$ is integrally closed precisely
when
\[
\surd{I}=IL\colon L^2,
\]
where $L=I \colon \surd{I}$.
We also show $($see Theorem~\ref{crit-2}$)$ that the equality
\[
I^2 \colon I =I
\]
gives a practical way of testing whether or not Gorenstein ideals of
codimension three are generically complete intersections.

We recall that an ideal $I$ satisfies a property {\it generically} if
$I_{\frak p}$ satisfies that property for every minimal prime ${\frak
p}$ of $I$.
We also say that an ideal $I$ in a local ring $R$ is a {\it complete
intersection of Goto-type} if $R$ is regular of dimension $d$ and $I=
(x_1, \ldots, x_{d-1}, x_d^t)$, where $x_1, \ldots, x_d$ is a regular
system of parameters $($see \cite[Theorem 1.1]{Goto}$)$.
The import of these special complete intersections is that they
are exactly the complete intersections which are integrally closed.

In Section 4, we then go on studying the following question: When $I$
integrally closed implies $I$ normal. Two cases are considered: the
integral closure of complete intersections of codimension $2$ in
Cohen--Macaulay local rings of arbitrary dimensions such that
$R_{\pp}$ has rational singularities in codimension two $($see
Theorem~\ref{codim2}$)$ and perfect Gorenstein  ideals of codimension
three in polynomials rings over fields of characteristic zero and with
certain restrictions on the $($local$)$ number of generators $($see
Theorem~\ref{normal}$)$. This makes the computational approach through
Rees algebras more accessible even though not inexpensive $($see
Remark~\ref{fromxcompu}, Example~\ref{computing}, and
Remark~\ref{improvement}$)$.

\section{Complete intersections}

As in the case of computation of radicals $($see \cite{EHV}$)$,
complete intersection ideals are the simplest case of study.
A careful  examination of \cite{Goto}, \cite{CPV}, and \cite{CP1} shows that
they deal with complementary aspects of integrality for such
ideals: While the first paper is concerned with integrally closed complete
intersections, the remaining two require that the ideal be {\em not}
integrally closed in their hypotheses.

We are going to frame these two  local aspects together into a global
integral closedness criterion that, when it fails, will still provide fresh
elements which are integral over the ideal.

\begin{theorem} \label{ccitest1}
Let $I$ be a height unmixed ideal in a Cohen--Macaulay ring $R$. Suppose
that $I$ is generically a complete intersection. Then $I$ is
integrally closed if and only if
\begin{equation}\label{rad}
\surd{I}=IL\colon L^2,
\end{equation}
where $L=I \colon \surd{I}$.
\end{theorem}
\demo
Let us assume first that $I$ is an integrally closed ideal. In order to
show that $($\ref{rad}$)$ holds it will suffice to prove it for any
$\p \in \mbox{Ass}(R/\surd{I})=\mbox{Ass}(R/I)$, given that
$\surd{I}\subseteq IL\colon L^2$, by the definition of $L$. If $\p
\in \mbox{Ass}(R/I)$, the strict inclusion $(\surd{I})_{\pp}=\p R_{\pp}
\subsetneq I_{\pp}L_{\pp} \colon L^2_{\pp}$ implies
$L^2_{\pp}=I_{\pp}L_{\pp}$. But $I_{\pp}$ is integrally closed, which forces
the equality $L_{\pp}=I_{\pp}$: a contradiction.

Conversely, suppose that $($\ref{rad}$)$ holds. Since $I$ is unmixed,
to show that $I=\overline{I}$ it suffices to check it at each
minimal prime $\pp$ of $I$.
First, if for such prime $R_{\pp}$ is not a regular local ring, the
ideal $L_{\pp}$ satisfies $I_{\pp}L_{\pp}= L^2_{\pp}$ by \cite{CP1},
which violates $($\ref{rad}$)$. Let $\p\in \mbox{Ass}(R/I)$, so that
$I_{\pp}=(f_1, f_2, \ldots, f_g)$, $g = \mbox{height}(I)$. If two
of the $f_i$'s belong to $\p^2R_{\pp}$, by
\cite[Theorem 2.1]{CPV}, one has that the ideal $L_{\pp}$ has
reduction number $1$ with respect to $I_{\pp}$, i.e.,
$L^2_{\pp}=I_{\pp}L_{\pp}$, which again would contradict $($\ref{rad}$)$.
Thus at most one of the $f_i$'s lies in the $\p^2R_{\pp}$. This means
that $I_{\pp}$ satisfies the conditions of \cite{Goto} and is
therefore integrally closed. \QED

\begin{remark}\label{methods}
To check the hypotheses on the ideal $I$ one can
proceed as follows $($see \cite{Eisenbudbook,xcompu}$)$.
Let $I$ be an ideal of codimension $m$.
\begin{description}[(a)]
\item[$($\mbox{\em a}$)$]
Suppose that $R$ is $($locally$)$ a Gorenstein ring and
let $J= (f_1, \ldots, f_m)$ be a subideal of $I$ of codimension $m$.
Then $I$ is height unmixed if and only if
\[
I= J \colon (J:I).
\]

\item[$($\mbox{\em b}$)$]
Let $R$ be a Cohen--Macaulay ring and let
\[
R^p \stackrel{\varphi}{\lar} R^q\lar I \rar 0
\]
be a presentation of $I$. Then $I$ is generically a complete intersection
if and only if the $($Fitting ideal$)$ ideal $I_{q-m}(\varphi)$ generated by
the minors of $\varphi$ of order $q-m$ has height at least $m+1$.
\end{description}
\end{remark}

\medskip

If $I$ is not integrally closed, elements that go into the  test of
Theorem~\ref{ccitest1} can be used to produce {\em new} elements in
its integral closure.  This occurs as follows

\begin{corollary}\label{ccitest2}
Let ${\p}_1, \ldots, {\p}_n$ be the minimal
prime ideals of $I$ listed in such a way that $I_{\pp}$ is
integrally closed for ${\p}= {\p}_i$ for $i\leq s$, but not
at the other primes. Set
\[
A= \surd{I}= {\p}_1 \cap \cdots \cap {\p}_n, \quad B= {\p}_1 \cap \cdots \cap
{\p}_s, \quad C= {\p}_{s+1} \cap \cdots \cap {\p}_n.
\]
Let $L=I \colon \surd{I}$. Then
\[
B=IL \colon L^2 \quad \mbox{and} \quad C = A \colon B.
\]
If $I$ is not integrally closed, that is if $B\neq A$, then
\[
H=I \colon C \neq  I \quad \mbox{and} \quad H^2=IH.
\]
\end{corollary}

\medskip

Note that one can arrange the computation in a manner that does not require
$\surd{I}$ directly. The following formulation mimics some of the
radical formulas of \cite{EHV}.

\begin{theorem}  \label{ccitest3}
Let $R$ be a ring of polynomials over a field of
characteristic zero and let $I$ be a height unmixed ideal. Let $J$
denote the Jacobian ideal of $I$. If $I$ is integrally closed then
\[
IJ \colon J=I.
\]
The converse holds if $I$ is generically a complete intersection.
\end{theorem}
\demo
We only consider the converse.
To show that $I$ is integrally closed it suffices to show that its
primary components are integrally closed. Localizing at minimal
primes of $I$ reduces the question to Theorem~\ref{ccitest1}.
But for any such prime ${\pp}$, $J_{\pp}$ is the generic socle of
$I_{\pp}$. \QED

\begin{example}\label{counter}
The converse of Theorem~\ref{ccitest3} does not hold without the
assumption of $I$ being generically a complete intersection.
For example, let $R$ be a regular local ring of dimension $2$ and maximal
ideal $\m=(x, y)$. Let $I$ be the Northcott ideal $(x^2, xy^4,
y^5)$. Its Jacobian ideal is $J=(xy^4, x^2y^3, y^8)$. It can be checked
that the condition $IJ \colon J = I$ is satisfied but $I$ is not
integrally closed. Indeed, the element
$xy^3 \not\in I$ belongs to the integral closure of $I$ as it satisfies
the monic equation $X^2 - y(x^2y^5) = 0$. It can actually be checked
that $\overline{I}=(x^2, xy^3, y^5)$.
\end{example}

\begin{remark}  When applying these methods where $L$ is either
$I\colon \surd{I}$, or the ordinary Jacobian ideal, the following
comment may be helpful.
If the test fails, that is if
\[
IL \colon L= I'\neq I,
\]
we could replace $I$ by $I'$ if the latter is height unmixed. Indeed
it has the same radical as $I$ so that its generic socle would be
\[
L' = I'\colon \surd{I},
\]
and we would test for
\[
I'L'\colon L' = I',
\]
and so on. There may be difficulties to this process. For instance, it was an
old question of Krull whether the integral closure of primary ideals
are still primary.  One of the authors $($see \cite{Huneke}$)$ gave a
counterexample in characteristic two but the characteristic zero case is
still open.
\end{remark}

\section{Gorenstein ideals}

%
%

The first result in this section describes, for a special class of ideals,
another practical way of telling when an ideal is generically a
complete intersection. More precisely, Theorem~\ref{crit-2}
shows that, in a Gorenstein ring $R$, a perfect Gorenstein ideal $I$ of
codimension three is generically a complete intersection if and only
if $I^2 \colon I =I$.
An immediate consequence of Theorem~\ref{ccitest1} and Theorem~\ref{crit-2}
is that for any such $I$ the following conditions are equivalent:
\begin{description}[(a)]
\item[$($\mbox{\em a}$)$]
$I$ is an integrally closed ideal;

\item[$($\mbox{\em b}$)$]
$I^2 \colon I = I$ and $\surd{I}=IL\colon L^2$, where $L=I \colon
\surd{I}$.
\end{description}
The rest of the section will be devoted to prove a generalization of
this fact. To be more precise, it will be shown in Theorem~\ref{new} that
for a generically Gorenstein ideal $I$ the condition
$\surd{I}=IL\colon L^2$, where $L=I \colon \surd{I}$, is necessary and
sufficient to guarantee that $I$ be integrally closed, without
restrictions on its height.
Some regularity assumptions on the ring $R$ are required as pointed
out in Example~\ref{regularity-assumpt}.

\begin{theorem}\label{crit-2}
Let $R$ be a Gorenstein ring and let $I$ be a perfect
Gorenstein ideal of codimension three. Then $I$ is generically a
complete intersection  if and only if
\[
I^2\colon I = I.
\]
\end{theorem}
\demo According to \cite{He}, for these ideals the conormal module
$I/I^2$ is Cohen--Macaulay and its associated primes are the
minimal primes of $I$. This means that we may localize at those
primes and reduce the question to the case of a local ring $R$ of
dimension $3$, from which it follows that if $I$ is a complete
intersection, then the asserted equality holds.

Conversely, suppose that $I^2\colon I = I$. This means that  $I/I^2$
is a faithful module of the Artinian, Gorenstein local ring $R/I$ and
therefore there is an embedding $R/I\hookrightarrow I/I^2$, which
leads to a decomposition
\[
I/I^2 \simeq R/I \oplus M,
\]
since $R/I$ is self-injective. We may assume that the image of $R/I$ in
$I/I^2$ has a lift $f$ which is a regular element of $I$. As in
\cite[Proof of Lemma 2]{W1}, this leads to a decomposition
\[
I/fI \simeq (f)/fI \oplus I/(f)\simeq R/I \oplus I/(f).
\]
This equality implies that as an $R/(f)$--module, $I/(f)$ is a
perfect ideal $($of codimension two$)$, whose first Betti number is the
same as the second Betti number of $I$. Thus $I/(f)$ is a perfect,
Gorenstein ideal of codimension two and therefore it is a complete
intersection, which means that $I$ is also a complete intersection. \QED

Corollary~\ref{Sunsook} generalizes a similar result contained in
\cite[Corollary 2.7]{SSV}, as the ideals under considerations are only
supposed to be integrally closed and not normal.

\begin{corollary}\label{Sunsook}
Let $R$ be a Gorenstein ring and let $I$ be a perfect
Gorenstein ideal of codimension three. If $I$ is integrally closed
then it is generically a complete intersection.
\end{corollary}
\demo
As $I$ is integrally closed one has, in particular, that $I^2 \colon I
= I$. The assertion then follows from Theorem~\ref{crit-2}.
\QED

\begin{remark}
A similar argument as the one in the proof of
Theorem~\ref{crit-2} will show that
if  $(R, \m)$ is a local ring of dimension $d\geq 3$ and $I$ is an
$\m$-primary ideal that is  perfect and Gorenstein then $I$ is a complete
intersection if and only if
\[
\bigwedge^{d-2}I/I^2
\]
is a faithful $R/I$--module.
\end{remark}

\begin{remark}\label{comparison}
From a computational point of view, the two methods $($the one involving
the Fitting ideals of $I$ and the one described in Theorem~\ref{crit-2}$)$
are essentially equivalent.
\end{remark}

%
%
%
%

The next two lemmas are crucial in the proof of Theorem~\ref{new}. The
first of them is about the existence of a {\em dual basis}.

\begin{lemma}\label{lem-dual}
Let $(R, \m)$ be a Noetherian local ring with embedding dimension $n$
at least two. Let $I$ be an $\m$-primary irreducible ideal contained
in $\m^2$. If $\m=(x_1, \ldots, x_n)$ then there exist $y_1,
\ldots, y_n$ such that for all $1 \leq i,j \leq n$
\begin{equation}\label{dualbasis}
x_iy_j \equiv \delta_{ij} \Delta \qquad \mbox{\rm mod\ } I,
\end{equation}
where $\Delta$ is a lift in $R$ of a socle generator of $R/I$ and
$\delta_{ij}$ denotes Kronecker's delta.
\end{lemma}
\demo
Note that $I \subset (x_1, \ldots, x_j^2, \ldots, x_n)$, as $I$ is
contained in ${\frak m}^2$.
For any $1 \leq j \leq n$ it will be enough to find $y_j
\in I \colon (x_1, \ldots, x_j^2, \ldots, x_n)$ such that $y_j
\not\in I \colon (x_1, \ldots, x_n)$. If not, one must have
$I \colon (x_1, \ldots, x_j^2, \ldots, x_n) = I \colon (x_1, \ldots,
x_n)$ and hence $(x_1, \ldots, x_j^2, \ldots, x_n) = (x_1, \ldots,
x_n)$, as $\mbox{\rm Hom}(\tratto, R/I)$ is a self-dualizing functor: a
contradiction.

The element $y_j$ has the property that $x_iy_j \in I$
for any $i \not= j$ and $x_jy_j \in I \colon (x_1,
\ldots, x_n)$. Hence we can write $y_jx_j = a_j+g_j \Delta$
with $a_j \in I$ and $g_j\not= 0$. However, $g_j \not\in
(x_1, \ldots, x_n)$ as otherwise $y_jx_j \in I$. But then $g_j$ is an
invertible element, so that $g_j^{-1}y_j$ will have all the required
properties as in $($\ref{dualbasis}$)$.
\QED

\begin{lemma}\label{new-0}
Let $(R, \m)$ be a Noetherian local ring with embedding dimension $n$
at least two. Suppose that $I$ is an $\m$-primary ideal
contained in $\m^2$ such that $R/I$ is Gorenstein. Letting $L=I \colon
\m$ then the following conditions hold:
\begin{description}[(a)]
\item[$(${\it a}$)$]
$L$ has reduction number one with respect to $I$, i.e., $L^2=IL$;

\item[$(${\it b}$)$]
$I \m = L \m$.
\end{description}
\end{lemma}
\demo
$(${\em a}$)$
As $L=(I, \Delta)$ one only needs to show that $\Delta^2 \in IL$. If
we let $\m=(x_1, \ldots, x_n)$, by Lemma~\ref{lem-dual} we can find
$y_1, \ldots, y_n$ and $a_1, \ldots, a_n \in I$ such that $\Delta =
x_iy_i +a_i$ for $1 \leq i \leq n$. As $n \geq 2$, we can write
\begin{eqnarray*}
\Delta^2 & = & (x_1y_1+a_1)(x_2y_2+a_2)=
                x_1y_1x_2y_2+x_1y_1a_2+a_1x_2y_2+a_1a_2 \\
         & = & (x_1y_2)(x_2y_1)+(x_1y_1)a_2+(x_2y_2)a_1+a_1a_2.
\end{eqnarray*}
Note that each term in the last sum belongs to the ideal
$I(I, \Delta)=I^2+I\Delta=IL$.

$(${\em b}$)$
We only need to show the inclusion $L \m \subseteq I \m$,
or better $\Delta \in I\m \colon \m$, as $L = (I, \Delta)$. But for
any $1 \leq i \leq n$ pick $j \not=i$ and write $\Delta = x_j
y_j + a_j$ with $a_j \in I$; hence
\[
x_i \Delta = x_i(x_j y_j + a_j) = x_j(x_iy_j) +x_ia_j \in \m I
\]
as desired.
\QED

\begin{theorem}\label{lemma-new}
Let $(R, {\frak m})$ be a Noetherian local ring and let $I$ be an ${\frak
m}$-primary ideal such that $R/I$ is Gorenstein. Then either
\begin{description}[(a)]
\item[$(${\it a}$)$]
there exists a minimal generating set $x_1, \ldots, x_n$ of the
maximal ideal ${\frak m}$ such that $I = (x_1, \ldots, x_{n-1},
x_n^s)$, or

\item[$(${\it b}$)$]
$L^2 = IL$, where $L = I \colon {\frak m}$.
\end{description}
\end{theorem}
\demo
Assume that $(${\it a}$)$ does not hold.
This
remains true after completing and further it suffices to prove $(${\it b}$)$
after completion. We can then write $R$ as a homomorphic image of a
regular local ring $S$ and it will suffice to prove $(${\it b}$)$ for
the pullback of $I$.
Notice that the pullback of $I$ to $S$ cannot satisfy $(${\it a}$)$.
Henceforth we may assume that $R$ is regular and
that $I$ does not satisfy $(${\it a}$)$. Let $t$ be the
vector space dimension of $(I+{\frak m}^2)/{\frak m}^2$ and
choose $x_1, \ldots, x_t$ in $I$ whose images span this vector space.
Replace $R$ by $R/(x_1, \ldots, x_t)$ and $I$ by $I/(x_1, \ldots,
x_t)$. We may then assume that $I$ is in ${\frak m}^2$. If the
dimension of $R$ is at least $2$ then Lemma~\ref{new-0} yields
$(${\it b}$)$. If not $R$ is either a field or a discrete valuation
ring and in this case $I$ satisfies $(${\it a}$)$, a contradiction.
\QED

\begin{remark}
Under the assumptions of Theorem~\ref{lemma-new} if we also assume that
$R$ is regular then we can restate $(${\it a}$)$ as
\begin{description}[(a')]
\item[$(${\it a$'$}$)$]
$I$ is a complete intersection of Goto-type.
\end{description}
\end{remark}

\begin{theorem} \label{new}
Let $I$ be a height unmixed ideal in a ring $R$. Suppose
that $I$ is a generically Gorenstein ideal and that $R_{\frak p}$ is a
regular local ring for all prime ideals ${\frak p}$ minimal over $I$.
Then $I$ is integrally closed if and only if
\begin{equation}\label{radnew}
\surd{I}=IL\colon L^2,
\end{equation}
where $L=I \colon \surd{I}$.
If either of these equivalent conditions hold, $I$ is generically a complete
intersection of Goto-type.
\end{theorem}
\demo
Let us assume first that $I$ is an integrally closed ideal. In order to
show that $($\ref{radnew}$)$ holds it will suffice to prove it for any
$\p \in \mbox{Ass}(R/\surd{I})=\mbox{Ass}(R/I)$, given that
$\surd{I}\subseteq IL\colon L^2$, by the definition of $L$. If $\p
\in \mbox{Ass}(R/I)$, the strict inclusion $(\surd{I})_{\pp}=\p R_{\pp}
\subsetneq I_{\pp}L_{\pp} \colon L^2_{\pp}$ implies
$L^2_{\pp}=I_{\pp}L_{\pp}$. But $I_{\pp}$ is integrally closed, which forces
the equality $L_{\pp}=I_{\pp}$: a contradiction.

Conversely, suppose that $($\ref{radnew}$)$ holds. Since $I$ is unmixed,
to show that $I=\overline{I}$ it suffices to check it at each
minimal prime $\pp$ of $I$. Without loss of generality we may assume
that $R$ is regular and $R/I$ is zero dimensional Gorenstein.
Then by Theorem~\ref{lemma-new} either $L^2=IL$ or $I$ is a complete
intersection of Goto-type. The first possibility cannot happen by
$($\ref{radnew}$)$ so $I$ must be a complete intersection of Goto-type
and is therefore integrally closed. \QED

\begin{example}\label{regularity-assumpt}
The regularity assumptions in Theorem~\ref{new} are required. For
example, let $S = k[X, Y, Z]$ be a polynomial ring in $3$ variables
over a field $k$. Let $R = S/(X^4+Y^4+Z^4)$ and let $x$, $y$, and $z$
denote the image of $X$, $Y$, and $Z$ in $R$, respectively. $R$ is
easily seen to be a two dimensional Gorenstein ring. Consider the
$R$-ideal $I = (x, y, z^2)$. Its lift $(X, Y, Z^2)$ back in $S$ is a
complete intersection of Goto-type, hence it is integrally
closed. However, $I$ itself is not an integrally closed $R$-ideal as
$L = (x, y, z)$ is integral over $I$. However, note that $(x, y, z)^2
\not= I(x, y,z)$.
\end{example}

\begin{corollary}\label{n=2}
With the same assumptions as in \mbox{\rm Theorem~\ref{new}}, if in
addition $I$ is an integrally closed ideal, in the linkage class of a complete
intersection with $R/I$ Gorenstein then
$I^2$ is integrally closed as well.
\end{corollary}
\demo
Since $R/I$ is Gorenstein and in the linkage class of a complete
intersection $R/I^2$ is Cohen--Macaulay, and in particular unmixed
$($see \cite{Buch-th}$)$. To prove $I^2$ is integrally closed it then
suffices to prove it locally at its minimal primes. After localizing
at any such prime $I$ is a complete intersection of Goto-type.
But then all powers of $I$ are also
generically complete intersections of Goto-type.
\QED

\begin{corollary}
Let $I$ be an unmixed ideal of a regular local ring $R$ and suppose that
$\overline{I}$ is Gorenstein. Then $I$ is integrally closed and it is
generically a complete intersection of Goto-type.
\end{corollary}
\demo
It is enough to show the equality $I=\overline{I}$ at the minimal associated
primes of $I$. Thus, after localizing, we may assume $R/\overline{I}$ to
be a zero dimensional Gorenstein ring. By Theorem~\ref{new} $\overline{I}$
is generically a complete intersection ideal of Goto-type.
In particular $\overline{I}$ cannot have proper reductions. As $I$ is
always a reduction of its integral closure, this forces the equality
$I=\overline{I}$, hence the claim.
\QED

\medskip

We conclude the section by pointing out that condition $($\ref{radnew}$)$ in
Theorem~\ref{new} is not sufficient to guarantee the integral closedness
in the case of ideals of type greater than or equal to $2$.

\begin{example}
Let $R$ be a regular local ring of dimension $2$ and maximal
ideal $\m=(x, y)$. The Northcott ideal $I = (x^2, xy^4, y^5)$ has
type $t=2$ and satisfies the condition $IL:L^2 = \m = \surd{I}$, with
$L = (x^2, xy^3, y^4)$. However $I$ is not integrally closed, as we
already observed in Example~\ref{counter}.
\end{example}

\section{Integral closedness and normality}

There are instances when it is more appropriate to prove that an ideal
is normal to avail ourselves of the Jacobian criterion.

\subsection*{Codimension two complete intersections}

The simplest case for tackling the question of computing the integral
closure of ideals
is certainly that of a codimension two complete intersection. Several
known facts come together to make the approach via Rees algebras amenable.
The methods to find the integral closure of affine domains become an
option $($see \cite{xcompu}$)$.

\begin{theorem}\label{codim2}
Let $R$ be a  Cohen--Macaulay local ring and suppose that $R_{\pp}$
has rational singularities for all codimension two prime ideals $\p$. If
$J=(a, b)$ is a complete intersection of codimension two, then
$\overline{J^n}= (\overline{J})^n$ for all $n\geq 1$.
\end{theorem}
\demo
For all $n \geq 1$ one has that
\[
(a, b)^{n-1} \overline{J} \subseteq (\overline{J})^n \subseteq
\overline{J^n},
\]
hence it will be enough to show that the equality
\begin{equation}\label{enough}
(a, b)^{n-1} \overline{J} = \overline{J^n}
\end{equation}
holds for all $n \geq 1$.

We claim that $R/(a, b)^{n-1} \overline{J}$ is unmixed of height
two. First of all, observe that
\[
(a, b)^{n-1}/(a, b)^{n-1}\overline{J} \simeq (a, b)^{n-1}/(a, b)^n
\otimes R/\overline{J} \simeq (R/(a, b))^n \otimes R/
\overline{J} \simeq (R /\overline{J})^n,
\]
as $(a, b)^{n-1}/(a, b)^n$ is isomorphic to the degree $n-1$ component of a
polynomial ring in two variables with coefficients in $R/(a, b)$.
On the other hand, one can consider the short exact sequence
\[
0 \rightarrow (a, b)^{n-1}/(a, b)^{n-1}\overline{J} \longrightarrow
R/(a, b)^{n-1}\overline{J} \longrightarrow R/(a, b)^{n-1} \rightarrow 0.
\]
The claim now follows from the fact that both $R/\overline{J}$ and
$R/(a, b)^{n-1}$ are unmixed $($for the first fact see \cite[Theorem
2.12]{Ratliff} and also \cite[Theorem 4.1]{McA}$)$.

Finally, it is enough to prove $($\ref{enough}$)$ at the associated
primes of $R/(a, b)^{n-1}\overline{J}$. In this case, however, the
result is true because of our assumption that $R_{\pp}$ has rational
singularities for all codimension two prime ideals $($see \cite{LT}$)$.
\QED

\begin{remark}\label{fromxcompu}
As mentioned earlier, an alternative approach to the computation of
the integral closure of an ideal $I$ is through the construction of
the integral closure of its Rees algebra $R[It]$. If $I=(a_1, \ldots,
a_n)$ then one can represent $R[It]$ as the quotient $R[T_1, \ldots,
T_n]/\QQ$, where $\QQ$ is the kernel of the map that sends $T_i$ to
$a_it$. If $R[It]$ is an affine domain over a field of characteristic
zero and $\mbox{Jac}$ denotes its Jacobian ideal, then the ring
\[
\mbox{Hom}_{R[It]}(\mbox{Jac}^{-1}, \mbox{Jac}^{-1})
\]
is guaranteed to be larger than $R[It]$ if the ring is not already
normal $($see \cite[Chapter 6]{xcompu}$)$. It also turns out that
\[
\mbox{Hom}_{R[It]}(\mbox{Jac}^{-1},
\mbox{Jac}^{-1})=\mbox{Hom}_{R[It]}((\mbox{Jac}^{-1})^{-1},
(\mbox{Jac}^{-1})^{-1})=(\mbox{Jac} \, \mbox{Jac}^{-1})^{-1},
\]
the latter being the inverse of the so called trace ideal of
$\mbox{Jac}$. Finally, to get the equations for $(\mbox{Jac} \,
\mbox{Jac}^{-1})^{-1}$ one then uses again \cite[Algorithm
6.2.1]{xcompu}.

Naturally, this process can be repeated several times until the
integral closure of $R[It]$ has been reached. The degree one component
of the final output gives the desired integral closure of $I$.
\end{remark}

\begin{example}\label{computing}
Let $k$ be a field of characteristic zero and let $J\subseteq R=k[x,
y]$ be the codimension two complete intersection
\[
J=(x^3+y^6, xy^3-y^5).
\]
Using the methods outlined in Remark~\ref{fromxcompu}, it turns out
that the integral closure of $J$ can be obtained in three steps.
After the first step has been completed one obtains the ideal
\[
J_1=(x^3+y^6, xy^3-y^5, y^8);
\]
after a second iteration one gets
\[
J_2=(x^3+y^6, xy^3-y^5, x^2y^2-y^6, y^7);
\]
finally, at the end of the last step, one has that the integral
closure of $J$ is given by the ideal
\[
J_3=\overline{J}=(xy^3-y^5, y^6, x^3, x^2y^2).
\]
Note that $\overline{J}$ is also a normal ideal.

For the records, despite the fact that the original setting for the
problem is a polynomial ring in $2$ variables over a field of
characteristic zero, overall we had to make use of additional 18
variables: quite a waste!
\end{example}

\begin{remark}\label{improvement}
In the case of the integral closure of complete intersections,
it has to be pointed out that some of the initial iterations of the
process described in Remark~\ref{fromxcompu} may be avoided by
using the results of \cite{CP2,PU}.

If the complete intersection $J$ is primary to the maximal ideal $\m$
of a Gorenstein ring $R$ and $J\subseteq \m^s$ but $J\not\subseteq
\m^{s+1}$ then one has an increasing sequence of ideals $I_k=J\colon
\m^k$ satisfying $I_k^2=JI_k$ for $k=1, \ldots, s$ if $\mbox{dim}(R)
\geq 3$ or for $k=1, \ldots s-1$ if $R$ is a regular local ring and
$\mbox{dim}(R)=2$. In particular this says that the $I_k$'s are
contained in the integral closure of $J$. Hence, instead of computing
the integral closure of $R[Jt]$ one may want to start directly from
$R[I_st]$ $($or $R[I_{s-1}t]$ if $R$ is a regular local ring and
$\mbox{dim}(R)=2)$.

If $J$ is not primary to the maximal ideal one may initially use
instead the sequence of ideals $I_k= J \colon (\surd{J})^{(k)} = J
\colon (\surd{J})^{k}$, provided that at each localization at the
associated primes of $J$ the conditions in \cite{CP2,PU} are satisfied
so that $I_k^2=JI_k$ for all $k$'s in the appropriate range.

In the specific case of the complete intersection $J$ of
Example~\ref{computing} one can easily compute $I_1=J \colon \m$ and
$I_2=J \colon \m^2$ and check that $J_1=I_1$ and $J_2=I_2$. Hence,
starting directly from a presentation of $R[I_2t]$ a
single iteration of the process as in Remark~\ref{fromxcompu} gives
the integral closure of $J$.
\end{remark}

\begin{remark}
In \cite{ES} D. Eisenbud and B. Sturmfels developed a nice theory about
binomial ideals. A natural question that one can raise on this matter asks:
Is the integral closure of a binomial ideal still a binomial
ideal? Despite what Example~\ref{computing} and similar other examples
would suggest, the answer is in general negative. We give below an example
in characteristic zero. We use the fact that an ideal of a polynomial ring
over a field is binomial if and only if some $($equivalently, every$)$
reduced Gr\"obner basis for the ideal consists of binomials $($see
\cite[Corollary 1.2]{ES}$)$.
\end{remark}

\begin{example}\label{bin-mon}
Let $R=k[x, y, z, w]$ be a ring with characteristic zero
and consider the ideal
\[
I=(x^2-xy, -xy+y^2, z^2-zw, -zw+w^2).
\]
It turns out that a primary decomposition for $I$ is given by
\[
I = (x-y, z-w) \cap (x-y, z^2, zw, w^2) \cap (z-w, x^2, xy, y^2) \cap
    (x, y^2, z, w^2).
\]
Hence the integral closure of $I$ is nothing but
\begin{eqnarray*}
\overline{I} & = & (x-y, z-w) \cap (x-y, z^2, zw, w^2) \cap (z-w, x^2, xy,
                   y^2) \\
             & = & (x^2-xy, -xy+y^2, z^2-zw, -zw+w^2, xz-yz-xw+yw).
\end{eqnarray*}
It can be checked that these generators for $\overline{I}$ are a
Gr\"obner basis and so it is not a binomial ideal.
The last fact can be easily checked if one chooses the order $x > z >
y >w$ and then shows that the $S$-resultant of each pair of the given
generating set of $\overline{I}$ reduces to zero modulo that same
generating set $($Buchberger criterion$)$.
\end{example}

\begin{remark}
It has to be remarked that with the following linear change of
variables
\[
X = x-y \qquad Y = y \qquad Z = z-w \qquad W = w
\]
the ideal $I$ of Example~\ref{bin-mon} becomes the following monomial ideal
\[
(X^2, -YX, Z^2, -WZ)
\]
so that its integral closure is again a monomial ideal, namely
\[
(X^2, -YX, Z^2, -WZ, XZ).
\]
\end{remark}

\subsection*{Gorenstein ideals}
We now address the normality for certain classes of
perfect Gorenstein ideals.

\begin{example}\label{example}
Let $k$ be a field of characteristic zero and $I \subset k[x, y,
z, w]$ the ideal generated by the Pfaffians of the five by five
skew-symmetric matrix
\[
\varphi=\left(
\begin{array}{ccccc}
0 & x & y & z & w \\
-x & 0 & x & y & z \\
-y & -x & 0 & x & y \\
-z & -y & -x & 0 & z \\
-w & -z & -y & -z & 0
\end{array}
\right).
\]
One has that
\[
I=(x^2-y^2+xz, xy-yz+xw,
xz-z^2+yw, xw, y^2-2xz)
\]
is a perfect Gorenstein ideal of codimension three. One can show that
\[
\surd I=(x-z, zw, yw, y^2-2z^2), \qquad L=(y, x,
z^2),
\]
and that the conditions specified in part $(${\em b}$)$ at the
beginning of Section~$3$ are satisfied. Thus, $I$ is an integrally closed
ideal. As a consequence of Corollary~\ref{coro} below, $I$ is
also a normal ideal.
%
\end{example}

\begin{theorem}\label{normal}
Let $R=k[x_1, \ldots, x_d]$ be a polynomial ring in $d$ variables over
a field $k$ of characteristic zero. Let $I \subset R$ be a Gorenstein
ideal defined by the Pfaffians of a $n\times n$ skew--symmetric matrix
$\varphi$ with linear entries. Suppose $n=d+1$ and that $I$ is a complete
intersection on the punctured spectrum. If $I$ is integrally closed it is
also normal.
\end{theorem}
\demo
$I$ is an ideal of analytic spread $\ell=d$, of reduction number $d-2$.
Its Rees algebra $R[It]$ is also Cohen--Macaulay and has a presentation
\[
R[It] = R[T_1, \ldots, T_n]/\L,
\]
with $\L=({\bf T}\cdot \varphi, f)$, where $f$ is the equation of
analytic dependence among the $n$ forms that generate $I$.
$($All of these facts can be traced to \cite{SUV}.$)$

To show that $I$ is normal, we estimate the codimension of the
non-normal locus of $R[It]$ and show that it is at least $2$. Since
$R[It]$ is Cohen--Macaulay this will suffice to prove that $I$ is normal.

Let $\N$ be the ideal that defines the non-normal locus of $R[It]$.
As $I$ is an integrally closed complete intersection on the punctured
spectrum of $R$, we may assume that the ideal $({\bf x})=(x_1, \ldots,
x_d)$ is contained in $\N$. We now look at the size of the Jacobian
ideal modulo $({\bf x})$. $\L$ is an ideal of height $n-1=d$, and the
corresponding $(n+d)\times (n+1)$ Jacobian matrix has the form

\medskip

\[
\left( \!
\begin{tabular}{c}
\parbox[t]{1.5in}{\parbox[t]{.89in}{
\fbox{$
\begin{array}{ccc}
& & \\
& B(\varphi) \\
& & \\
\end{array} $}
 }\parbox[t]{.4in}{$
\begin{array}{ccc}
\ &    0   & \ \\
\ & \vdots & \ \\
\ &    0   & \
\end{array}
$} \\
\parbox[t]{.89in}{\mbox{$
\begin{array}{cccc}
& &           & \\
& & \ \varphi & \\
& &           &
\end{array}
$}}\parbox[t]{.4in}{
\fbox{$
\begin{array}{c}
\displaystyle\frac{\partial f}{\partial T_1}^{\ } \\
\vdots \\
\displaystyle\frac{\partial f}{\partial T_n}_{\ } \\
\end{array} $}
} }
\end{tabular}
\!\! \right)
\]

\medskip

\noindent
where $B(\varphi)$ is the $($Jacobian dual$)$ matrix of the linear forms
in the $T_i$'s such that ${\bf T}\cdot \varphi = {\bf x}\cdot B(\varphi)$.
Note that $B(\varphi)$ is a $d\times n$ matrix.

The Jacobian ideal contains the product of the minors of order $d-1$ of
$B(\varphi)$ by the ideal generated by the partial derivatives of $f$.
Since $f$ is an irreducible polynomial and $k$ has characteristic
zero,  the latter has height at least
$2$. On the other hand, based on an argument from \cite[Proposition
2.4]{Morey}, one has that $\mbox{\rm height } I_{d-1}(B(\varphi))\geq
2$. Hence we have an ideal of forms in the variables $T_i$'s to add to
the ideal $({\bf x})$; both are contained in $\N$ and from the
previous calculations it follows that
\[
\mbox{\rm height } \N\geq d+2= \mbox{\rm height } \L + 2,
\]
as claimed.  \QED

\begin{corollary}\label{coro}
Let $R=k[x_1, x_2, x_3, x_4]$ be a polynomial ring in $4$ variables over
a field $k$ of characteristic zero. Let $I \subset R$ be a Gorenstein
ideal defined by the Pfaffians of a five by five skew--symmetric matrix
$\varphi$ with linear entries. If $I$ is integrally closed it is
also normal.
\end{corollary}
\demo
By Theorem~\ref{crit-2}, $I$ is a complete intersection on the
punctured spectrum, hence Theorem~\ref{normal} applies. \QED

\begin{example}
The hypothesis that the entries of the matrix
$\varphi$ are linear forms is necessary.
Let $k$ be a field of characteristic zero and $I \subset k[x, y,
z, w]$ be the codimension three Gorenstein ideal $I$ generated by the
Pfaffians of the five by five skew-symmetric matrix $($see \cite[Example
8.3.3]{xcompu}$)$
\[
\varphi =
\left(
\begin{array}{ccccc}
0   & -x^2 & -y^2 & -z^2 & -w^2 \\
x^2 & 0    & -w^2 & -xy  & -z^2 \\
y^2 & w^2  & 0    & -x^2 & -xy  \\
z^2 & xy   & x^2  & 0    & -y^2 \\
w^2 & z^2  & xy   & y^2  & 0
\end{array}
\right),
\]
where $x, y, z, w$ are variables. It can be checked that $IL \colon
L^2 = \surd{I}$ where $L=I \colon \surd{I}$. Hence $I$ is integrally closed.
On the other hand, this ideal is not normal, a fact that can be shown
as follows. Since $I$ is generated by five forms $f_1, \ldots, f_5$ of
the same degree, the Rees algebra of $I$ decomposes as
\[
R[It]=k[f_1, \ldots, f_5] \oplus \L.
\]
To prove that $R[It]$ is not normal it suffices to show that $k[f_1,
\ldots, f_5]$ is not normal. But $k[f_1,\ldots, f_5]$ can be easily
seen to be an hypersurface ring and an application of the Jacobian
criterion shows that it is not normal.
\end{example}

\end{document}